\documentclass[amscd,amssymb,11pt,reqno]{amsart}

\newtheorem{Thm}{Theorem}[section]

\newtheorem{Cor}[Thm]{Corollary}
\newtheorem{Lem}[Thm]{Lemma}
\newtheorem{Prop}[Thm]{Proposition}

\newtheorem{Rmk}[Thm]{Remark}

\usepackage{graphicx}
\usepackage{latexsym}
\begin{document}

%\begin{Large}

\vspace{1.5 cm}

\title[The Busemann-Petty problem  in hyperbolic and spherical spaces]
      {The Busemann-Petty problem\\in hyperbolic and spherical spaces}

\author{  V.Yaskin}

\address{Department of Mathematics, University of Missouri, Columbia, MO 65211, USA.}

\email{yaskinv@math.missouri.edu}

%\begin{center}{\bf Sections} \end{center}

\begin{abstract}

The Busemann-Petty problem asks whether origin-symmetric convex
bodies in $\mathbb{R}^n$ with smaller central hyperplane sections
necessarily have smaller $n$-dimensional volume. It is known that
the answer to this problem is affirmative if $n\le 4$ and negative
if $n\ge 5$. We study this problem in hyperbolic and spherical
spaces.

\end{abstract}

\maketitle

\section{Introduction}

The Busemann-Petty problem asks the following question. Given two
convex origin-symmetric bodies $K$ and $L$ in $\mathbb{R}^n$ such
that
$$ \mathrm{vol}_{n-1}(K\cap H)\le \mathrm{vol}_{n-1} (L\cap H)$$
for every central hyperplane $H$ in $\mathbb{R}^n$, does it follow
that
$$\mathrm{vol}_n(K)\le \mathrm{vol}_n (L)?$$

The answer to this problem in $\mathbb{R}^n$ is known to be affirmative if $n\le 4$ and negative if
$n\ge 5$. The solution appeared as the result of work of many mathematicians (see \cite{GKS} or
\cite{Zh}  for  historical details).

In this paper we consider the Busemann-Petty problem in hyperbolic and spherical spaces in place of
the Euclidean space. We prove

\begin{Thm}
Let $K$ and $L$ be centrally symmetric  convex bodies in the spherical space $\mathbb{S}^n$, $n\le
4$ (more precisely in a hemisphere) such that
\begin{eqnarray}\label{eqn:condition}
 \mathrm{vol}_{n-1}(K\cap H)\le \mathrm{vol}_{n-1} (L\cap H)
\end{eqnarray}
for every central totally-geodesic hyperplane $H$ in $ \mathbb{S}^n$. Then $$\mathrm{vol}_n(K)\le
\mathrm{vol}_n (L).$$ On the other hand, if $n\ge 5$ there are convex symmetric bodies $K$,
$L\subset \mathbb{S}^n$ that satisfy {\rm (\ref{eqn:condition})} but $
\mathrm{vol}_n(K)>\mathrm{vol}_n(L).$
\end{Thm}

So, the answer to the Busemann-Petty in $\mathbb{S}^n$ is exactly the same as in the Euclidean
space. However, the situation in  the hyperbolic space is  different. Trivially, the answer is
affirmative if $n=2$, since the condition (\ref{eqn:condition}) in this case is equivalent to $K
\subseteq L$, but for  higher dimensions we have the following:
\begin{Thm}
There are convex centrally symmetric bodies $K$, $L\subset \mathbb{H}^n$, $n\ge3$  that satisfy the
condition
\begin{eqnarray*}
 \mathrm{vol}_{n-1}(K\cap H)\le \mathrm{vol}_{n-1} (L\cap H)
\end{eqnarray*} for every central totally-geodesic hyperplane $H$ in $\mathbb{H}^n$,
but $\mathrm{vol}_n(K)> \mathrm{vol}_n(L).$
\end{Thm}

The idea to find analogs of known results in non-Euclidean spaces is not new. For example in
\cite{GHS} the authors study intrinsic volumes in hyperbolic and spherical spaces. The
Brunn-Minkowski inequality in different spaces is discussed in \cite{G}. Also a number of papers is
concerned with other generalizations of the Busemann-Petty problem. In our proof we will be using
results from \cite{Zv}, where A.Zvavitch studied the Busemann-Petty problem for arbitrary measures.
For other generalizations of the Busemann-Petty problem  see \cite{BZ}, \cite{K3}, \cite{K4},
\cite{K5}, \cite{K6}, \cite{RZ}, \cite{KYY}.

\section{Preliminaries}
Let $\mathbb{S}^n$ be the unit sphere  in $\mathbb{R}^{n+1}$. Using the stereographic projection
(from the north pole onto the hyperplane containing the equator) we can think of it as
$\mathbb{R}^{n}$ equipped with the metric of constant curvature $+1$:
$$ds^2=4 \frac{dx_1^2+\cdots +dx_n^2}{(1+(x_1^2+\cdots+x_n^2))^{2}},$$
where $x_1$,..., $x_n$ are the standard Euclidean coordinates in $\mathbb{R}^{n}$. (See \cite[\S 9,
\S 10]{DFN}, and \cite[\S 4.5]{R}  for details about the spherical and hyperbolic spaces). It is
well-known that geodesic lines on the sphere are great circles. Later on, in order to define
convexity, we will need the uniqueness property of geodesics joining given 2 points. But this is
not the case on the sphere. However if we restrict ourselves to an open hemisphere, then for any
two points there exists a unique geodesic segment connecting them. Under the stereographic
projection the open south hemisphere gets mapped onto the open unit ball $B^n$ in $\mathbb{R}^{n}$.
This is the model we will be working in. The geodesics in this model are arcs of the circles
intersecting the boundary of the ball $B^n$ %$x_1^2+\cdots + x_n^2=1$
in antipodal points and straight lines through the origin.

Also it is well-known that the hyperbolic space $\mathbb{H}^n$ can be identified with the interior
of the unit ball in $\mathbb{R}^{n}$ with the metric:

$$ds^2=4 \frac{dx_1^2+\cdots +dx_n^2}{(1-(x_1^2+\cdots+x_n^2))^{2}}.$$
This is the Poincar\'{e} model of the hyperbolic space in the ball.  Note that it can be also
obtained from the pseudeosphere in the Lorentzian space via the stereographic projection. The
geodesic lines in this model are arcs of the circles orthogonal to the boundary of the ball $B^n$ %sphere $x_1^2+\cdots +x_n^2=1$
and straight lines through the origin.

Since both geometries are defined in the unit ball in $\mathbb{R}^{n}$, we will treat them
simultaneously, considering the open ball $B^n\subset \mathbb{R}^n$ with the metric
\begin{eqnarray}\label{eqn:metric}
ds^2=4 \frac{dx_1^2+\cdots +dx_n^2}{(1+\delta\ (x_1^2+\cdots+x_n^2))^{2}},
\end{eqnarray}
where $\delta = -1$ for the hyperbolic case, $+1$ for the spherical space. In addition if we
consider $\delta=0$ we get the original case of the Euclidean space.

The  definition of convexity in hyperbolic and spherical spaces (recall that we work in an open
hemisphere) is analogous to that in the Euclidean space (see \cite[Chapter I, \S 12]{P}).  A body
$K$ (compact set with non-empty interior) is called {\bf convex} if for every pair of points in $K$
the geodesic segment joining them also belongs to the body $K$. For our definition of convexity in
$\mathbb{S}^n$ it is crucial that we work in an open hemisphere, since in this case we have a
unique geodesic segment through any two points.

Let $K$ be a body in the open unit ball $B^n$. In order to distinguish between different types of
convexity we will adopt the following system of notations.  The body $K$ is called s-convex (or
$+1$-convex), if it is convex in the spherical metric defined in the ball $B^n$. Similarly it is
called h-convex (or $-1$-convex) if it is convex with respect to the hyperbolic metric.  e-convex
bodies (or 0-convex) are the bodies convex in the usual Euclidean sense. Analogously
s-(h-,e-)geodesics are the straight lines of the spherical (hyperbolic, Euclidean) metric. (In this
terminology we follow \cite{MP}. Note that in the literature there are other definitions of
h-convexity or $\delta$-convexity which have absolutely different meaning).

Shown below  are some examples of convex hulls of 4 points with respect to hyperbolic, Euclidean
and spherical metrics correspondingly.

\begin{center}
\includegraphics{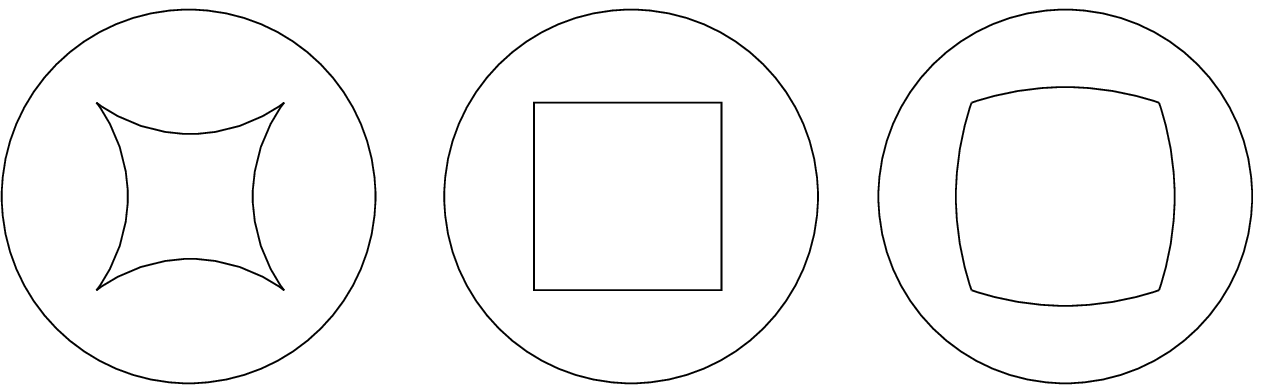}
\vspace{0.2cm}

Figure 1

\vspace{0.2cm}

\end{center}

Clearly, any s-convex body containing the origin is also e-convex and any e-convex body containing
the origin is h-convex. (See for example \cite{MP}).

A submanifold $\mathcal F$ in a Riemannian space $\mathcal R$ is called {\bf totally geodesic} if
every geodesic in $\mathcal F$ is also a geodesic in the space $\mathcal R$. In the Euclidean space
the totally geodesic submanifolds are Euclidean planes, on the sphere they are great subspheres. In
the  Poincar\'{e} model of the hyperbolic space described above the totally geodesic submanifolds
are represented by the spheres orthogonal to the boundary of the unit ball $B^n$ and Euclidean
planes through the origin. In a sense, totally geodesic submanifolds are analogs of Euclidean
planes in Riemannian spaces. For elementary properties of totally geodesic submanifolds see
\cite[Chap.5, \S 5]{A}.

The {\bf Minkowski functional} of a star-shaped origin-symmetric body $K\subset \mathbb R^n$ is
defined as
$$||x||_K=\min \{a\ge 0: x \in aK \}.$$
The {\bf radial function} of $K$ is given by $\rho_K(x) =||x||_K^{-1}$. If $x\in S^{n-1}$ then the
radial function $\rho_K(x)$ is the Euclidean distance from the origin to the boundary of $K$ in the
direction of $x$.

For a  centrally-symmetric $\delta$-convex body $K\in B^n$ ($\delta=0,1,-1$) consider the section
of $K$ by the hypersurface $\xi^\perp =\{\langle x,\xi\rangle=0\}$, where $\langle
\cdot,\cdot\rangle$ is the Euclidean scalar product. Clearly such a hypersurface  is a totally
geodesic hyperplane in the metric (\ref{eqn:metric}) for any $\delta=0,1,-1$. This hyperplane
passes through the origin with the normal vector $\xi$.

The volume element of the metric (\ref{eqn:metric}) equals

$$d\mu_n=2^n \frac{dx_1\cdots dx_n}{(1+\delta\ (x_1^2+\cdots+x_n^2))^{n}}=2^{n} \frac{dx}{(1+\delta\ |x|^2)^{n}}.$$
Therefore the volume of a body $K$ is given by the formula:
$$\mathrm{vol}_n(K)= \int_{K} d\mu_{n}=2^n \int_K \frac{dx}{(1+\delta\ |x|^2)^{n}}.$$
Note that in polar coordinates the latter formula looks as follows:
\begin{eqnarray}\label{eqn:polarvolume}
\mathrm{vol}_n(K)=2^n\int_{S^{n-1}}\int_{0}^{||\theta||^{-1}_K} \frac{r^{n-1}}{(1+\delta\
r^2)^{n}}dr\ d\theta.
\end{eqnarray}

 Similarly the volume element of the hypersurface $\xi^\perp$ is
$$d\mu_{n-1}=2^{n-1} \frac{dx}{(1+\delta\ |x|^2)^{n-1}},$$
therefore the $(n-1)$-volume of the section of $K$ by the hyperplane $\xi^\perp$ is given by the
formula:
$$S_{K}(\xi)=\int_{K\cap\langle x,\xi\rangle=0} d\mu_{n-1}=2^{n-1} \int_{K\cap\langle x,\xi\rangle=0} \frac{dx}{(1+\delta\ |x|^2)^{n-1}}.$$

One of the tools of this paper is the Fourier transform of distributions. The Fourier transform of
a distribution $f$ is defined by $\langle\hat{f}, {\phi}\rangle= \langle f, \hat{\phi} \rangle$ for
every test function $\phi$ from the space $ \mathcal{S}$ of rapidly decreasing infinitely
differentiable functions on $\mathbb R^n$.

A distribution is called {\bf positive definite} if for every test function $\phi$
$$\langle f, \phi \ast \overline{\phi(-x)}\rangle \ge 0.$$
By L.Schwartz's generalization of Bochner's theorem, a distribution is positive definite if and
only if its Fourier transform is a positive distribution (in the sense that $\langle \hat{f},\phi
\rangle \ge 0$ for every non-negative test function $\phi$; see, for example, \cite[p.152]{GV}).

The spherical Radon transform $R:C(S^{n-1})\to C(S^{n-1})$ is defined by

$$Rf(\xi)=\int_{S^{n-1}\cap \xi^\perp} f(x) dx.$$

The following Lemma, due to A.Koldobsky, gives a relation between the spherical Radon transform and
the Fourier transform.

\begin{Lem}\label{Lem:Radon}{\rm (\cite{K1}, Lemma 4)}
Let $g(x)$ be an even homogeneous function of degree $-n+1$ on $\mathbb{R}^n\setminus \{0\}$,
$n>1$, so that $g(x)|_{S^{n-1}}\in C({S^{n-1}})$ then

$$R g(\xi)=\frac{1}{\pi} \hat{g}(\xi), \qquad \forall \xi \in S^{n-1}.$$
\end{Lem}

The latter equality means that $\hat{g}$ is a homogeneous function of degree $-1$ on
$\mathbb{R}^n$, whose values on $S^{n-1}$ are equal to $Rg$.

 Now we derive a formula for the function $S_{K}(\xi)$ using the Fourier transform, similar to
\cite{Zv}. For $\delta=0$ this is the formula from \cite[Theorem 1]{K1}.

\begin{Lem}\label{Lem:Main}
 Let $K$ be an origin-symmetric $\delta$-convex body in $B^n$ with
 Minkowski functional $||\cdot||_K$. Let $\xi\in S^{n-1}$ and $\xi^\perp$ be the hyperplane through the origin orthogonal to $\xi$.
 Then the volume of the section of the body $K$ by the hyperplane $\xi^\perp$
 in the metric (\ref{eqn:metric}) equals
$$S_{K}(\xi)=\frac{2^{n-1}}{\pi}\left(|x|^{-n+1}_2
\int_0^{\frac{|x|}{||x||_K}}\frac{r^{n-2}}{(1+\delta\ r^2)^{n-1}}\ dr \right)^\wedge(\xi).$$

\end{Lem}
\noindent{\bf Proof.} Passing to spherical coordinates we get:
\begin{eqnarray*}
S_{K}(\xi)  &=&2^{n-1} \int_{\xi^\perp} \chi(||x||_K) \frac{dx}{(1+\delta\ |x|^2)^{n-1}}=\\
            &=&2^{n-1} \int_{S^{n-1}\cap\xi^\perp}
            \int_0^{||\theta||_K^{-1}} \frac{r^{n-2} dr}{(1+\delta\ r^2)^{n-1}}\ d\theta.
\end{eqnarray*}

We can  rewrite the integral above as follows (note that $|x|=1$, since $x\in S^{n-1}$):
\begin{eqnarray*}
S_{K}(\xi)  &=&2^{n-1} \int_{S^{n-1}\cap\xi^\perp} |x|^{-n+1}
            \int_0^{|x|/||x||_K} \frac{r^{n-2} dr}{(1+\delta\ r^2)^{n-1}}\ dx.
\end{eqnarray*}
The function  under the spherical integral is a homogeneous function of $x$ of degree $-n+1$ and
therefore by Lemma \ref{Lem:Radon}:

$$S_{K}(\xi)=\frac{2^{n-1}}{\pi}\left(|x|^{-n+1}_2
\int_0^{\frac{|x|}{||x||_K}}\frac{r^{n-2}}{(1+\delta\ r^2)^{n-1}}\ dr \right)^\wedge(\xi).$$

\qed

\section{Proofs of main results}

First we construct counterexamples to the Busemann-Petty problem in $\mathbb H^n$ and $\mathbb S^n$
for $n\ge 5$.

\begin{Thm}\label{Thm:negative}
There exist convex origin-symmetric bodies $K$ and $L$ in $\mathbb S^n$ (or $\mathbb H^n$), $n\ge
5$ such that
\begin{eqnarray*}
 \mathrm{vol}_{n-1}(K\cap H)\le \mathrm{vol}_{n-1} (L\cap H)
\end{eqnarray*}
for every central hyperplane, but $\mathrm{vol}_n(K)> \mathrm{vol}_n(L).$
\end{Thm}
\noindent{\bf Proof.} We will show the proof only for the case of the spherical space, the
hyperbolic case is similar. The idea here is to use the property that any Riemannian space locally
looks as ``almost" Euclidean.

Let $K$ and $L$ be  convex origin-symmetric bodies in $\mathbb R^n$ that give a counterexample to
the original Busemann-Petty problem. That is
\begin{equation}\label{eqn:sections}
\mathrm{EVol}_{n-1}(K\cap H)\le \mathrm{EVol}_{n-1} (L\cap H)
\end{equation} for every central
hyperplane $H$, but
\begin{equation}\label{eqn:vol}\mathrm{EVol}_n(L)< \mathrm{EVol}_n(K).
\end{equation} (Here we denote the usual Euclidean volume by EVol to avoid confusion with the spherical
volume.)

In fact,  since the inequality (\ref{eqn:vol}) is strict, we can dilate one of the bodies a little
to make the inequality (\ref{eqn:sections}) strict. Recall also, that in the original
counterexample the body $L$ was strictly convex, and the body $K$ was obtained from the body $L$ by
small perturbations. Note that  $K$ can also be made strictly convex.

In view of the latter remarks,  we will assume that $K$ and $L$ are  strictly convex
origin-symmetric
bodies that satisfy the strict version of (\ref{eqn:sections}). Moreover, %since we have  strict inequalities for volumes,
there exists an $\epsilon>0$ such that
$$\mathrm{EVol}_{n-1}(K\cap H)< (1-\epsilon)\mathrm{EVol}_{n-1} (L\cap H)$$ for all $H$ and
$$\mathrm{EVol}_n(L)< (1-\epsilon) \mathrm{EVol}_n(K).$$ Clearly, any
dilations $\alpha K$ and $\alpha L$ also provide a counterexample. We can take $\alpha$ so small
that both bodies $K$ and $L$ lie in a ball of radius $r$ that satisfies the inequality:
$$1-\epsilon \le \frac{1}{(1+r^2)^n}\le 1.$$

Now the volumes of the bodies $K$ and $L$  in the  spherical metric are related by the inequality:
\begin{eqnarray*}
\mathrm{vol}_n(L)   &=&2^{n}\int_{L}\frac{dx}{(1+|x|^2)^n}\le
                     2^{n}\int_{L} dx=
                    2^n \mathrm{EVol}_n(L)< \\
                    &<&(1-\epsilon) 2^n\mathrm{EVol}_n(K)=
                    (1-\epsilon) \ 2^{n}\int_{K} dx\le\\
                    &\le&2^{n}\int_{K}\frac{dx}{(1+|x|^2)^n}=\mathrm{vol}_n(K).
\end{eqnarray*}

Analogously, for the volumes of sections we have
\begin{eqnarray*}
\mathrm{vol}_{n-1}(K\cap \xi^\perp) &=&2^{n-1} \int_{K\cap\langle x,\xi\rangle=0} \frac{dx}{(1+\ |x|^2)^{n-1}}\le\\
                                    &\le& 2^{n-1} \int_{K\cap\langle x,\xi\rangle=0}dx <\\
                                    &<&(1-\epsilon) 2^{n-1} \int_{L\cap\langle x,\xi\rangle=0}dx\le\\
                                    &\le& 2^{n-1} \int_{L\cap\langle
x,\xi\rangle=0} \frac{dx}{(1+\ |x|^2)^{n-1}}= \mathrm{vol}_{n-1} (L\cap \xi^\perp).
\end{eqnarray*}

To finish the proof we only need to show that if $K$ is a strictly e-convex body, then $\alpha K$
is s-convex for sufficiently small $\alpha$. Consider the boundary of the body $K$. Define
$$k=\min\{k_i(x): x\in \partial K, \ i=1,...,n-1\},$$
where $k_i(x)$, $i=1$,..., $n-1$, are the principal curvatures at the point $x$ on the boundary of
$K$. Since $K$ is strictly e-convex the quantity defined above is strictly positive: $k>0$. For the
body $\alpha K$ it is equal to $\displaystyle k/\alpha$. On the other hand in a small neighborhood
of the origin the totally geodesic s-planes are the spheres with almost zero curvature (from the
Euclidean point of view). Consider all the spheres, which are totally geodesic in the spherical
metric and tangent to the body $\alpha K$, and let $R$ be the smallest radius of all such spheres.
We can choose an $\alpha$ so small that
$$k/\alpha> 1/R$$ and therefore the body $\alpha K$ lies on one side with
respect to any tangent totally geodesic s-hyperplane. Hence $\alpha K$ is s-convex.

The situation in the hyperbolic space is even easier since every e-convex body containing the
origin is also h-convex.

\qed

In 1988 E.Lutwak \cite{L} introduced the concept of intersection body and proved that the
Busemann-Petty problem has affirmative answer if the body with smaller sections is an intersection
body. Later, in \cite{K2} A.Koldobsky proved that a body $K$ is an intersection body if and only if
$||x||_K^{-1}$ is a positive definite distribution. Then in \cite{K3} A.Koldobsky generalized
Lutwak's connection using the following Parseval's formula on the sphere:

\begin{Lem}\label{Lem:Parseval}
If $K$ and $L$ are origin symmetric infinitely smooth bodies in $\mathbb{R}^n$ and $0<p<n$,
then

$$\int_{S^{n-1}} \left(||x||_K^{-p}\right)^\wedge (\xi) \left(||x||_L^{-n+p}\right)^\wedge
(\xi)d\xi= (2\pi)^n \int_{S^{n-1}} ||x||_K^{-p}||x||_L^{-n+p} dx.$$
\end{Lem}

 In fact we will be using the following version of this Lemma, see \cite[Corollary 1]{K3}.
\begin{Cor}\label{Cor:Parseval}
Let $f$ and $g$ be functions on $\mathbb{R}^n$, continuous on $S^{n-1}$ and homogeneous of degree
$-1$ and $-n+1$ respectively. Suppose that $f$ represents a positive definite distribution. Then
there exists a measure $\gamma_0$ on $S^{n-1}$ such that
\begin{equation*}
\int_{S^{n-1}} \widehat{g} (\theta)\, d\gamma_0(\theta)= (2\pi)^n \int_{S^{n-1}} f(\theta)\,
g(\theta)\, d\theta.
\end{equation*}
\end{Cor}

 Later, A.Zvavitch (\cite{Zv}) solved
the Busemann-Petty problem for arbitrary measures. Namely, let $f_n(x)$ be a locally integrable
function on $\mathbb{R}^n$, and $f_{n-1}(x)$ a function on $\mathbb{R}^n$, locally integrable on
central hyperplanes. Then  let $\mu_n$ be the measure on $\mathbb{R}^n$ with density $f_n(x)$ and
$\mu_{n-1}$ be the $(n-1)$-dimensional measure on central hyperplanes with density $f_{n-1}(x)$
such that $t\frac{f_n(tx)}{f_{n-1}(tx)}$ is an increasing function of $t$ for any fixed $x$. Then
if
$$||x||_K^{-1}\frac{f_n(\frac{x}{||x||_K})}{f_{n-1}(\frac{x}{||x||_K})}$$
is a positive definite distribution on $\mathbb{R}^n$ then the Busemann-Petty problem for these
measures has affirmative answer, i.e. $\mu_{n-1}(K\cap\xi^\perp)\le \mu_{n-1}(L\cap\xi^\perp)$
implies $\mu_n(K)\le \mu_n(L)$. Our next result is a particular case of  Zvavitch's theorem, but
for the sake of completeness we include  a proof.
\begin{Thm}\label{Thm:Main}

Let $K$ and $L$ be $\delta-$convex origin-symmetric  bodies in $B^n$ such that $\displaystyle
\frac{||x||^{-1}_K}{1+\delta\ (\frac{|x|}{||x||_K})^2}$ is a positive definite distribution. If
$$\mathrm{vol}_{n-1}(K\cap H)\le \mathrm{vol}_{n-1}(L\cap H)$$ for every totally geodesic
hyperplane through the origin, then
$$\mathrm{vol}_{n}(K)\le \mathrm{vol}_{n}(L).$$

\end{Thm}

\noindent{\bf Proof.} Let us first prove the following elementary inequality (cf. Zvavitch,
\cite{Zv}). For any $a,b\in (0,1)$
\begin{eqnarray*}
\frac{a}{1+\delta\ a^2}\int_a^b \frac{r^{n-2}}{(1+\delta\ r^2)^{n-1}}dr\le \int_a^b
\frac{r^{n-1}}{(1+\delta\ r^2)^{n}}dr.
\end{eqnarray*}

Indeed, since the function $\displaystyle \frac{r}{1+\delta\ r^2}$ is increasing on the interval
$(0,1)$ we have the following
\begin{eqnarray*}
\frac{a}{1+\delta\ a^2}\int_a^b \frac{r^{n-2}}{(1+\delta\ r^2)^{n-1}}dr
&=&\int_a^b\frac{r^{n-1}}{(1+\delta\ r^2)^{n}}\frac{a}{1+\delta\ a^2}\left(\frac{r}{1+\delta\
r^2}\right)^{-1} dr \\
&\le& \int_a^b \frac{r^{n-1}}{(1+\delta\ r^2)^{n}}dr.
\end{eqnarray*}
Note that latter inequality does not require that $a\le b$.

 Using the previous inequality with $a=||x||^{-1}_K$ and $b=||x||^{-1}_L$ we get
\begin{eqnarray*}
\int_{S^{n-1}}\frac{||x||^{-1}_K}{1+\delta\ ||x||^{-2}_K}\int_{||x||^{-1}_K}^{||x||^{-1}_L}
\frac{r^{n-2}}{(1+\delta\ r^2)^{n-1}}drdx\le \int_{S^{n-1}}\int_{||x||^{-1}_K}^{||x||^{-1}_L}
\frac{r^{n-1}}{(1+\delta\ r^2)^{n}}drdx.
\end{eqnarray*}

Suppose we can show that the left-hand side is non-negative, then it will follow that
\begin{eqnarray*}
\int_{S^{n-1}}\int_{0}^{||x||^{-1}_K} \frac{r^{n-1}}{(1+\delta\
r^2)^{n}}drdx\le\int_{S^{n-1}}\int_{0}^{||x||^{-1}_L} \frac{r^{n-1}}{(1+\delta\ r^2)^{n}}drdx,
\end{eqnarray*}
that is $\mbox{vol}_{n}(K)\le \mbox{vol}_{n}(L),$ see the polar formula (\ref{eqn:polarvolume}).

So we only need to show that
\begin{eqnarray*}
\int_{S^{n-1}}\frac{||x||^{-1}_K}{1+\delta\ ||x||^{-2}_K}\int_{0}^{||x||^{-1}_K}
\frac{r^{n-2}}{(1+\delta\ r^2)^{n-1}}drdx\le \hspace{4cm}\\ \hspace{4cm}\le
\int_{S^{n-1}}\frac{||x||^{-1}_K}{1+\delta\ ||x||^{-2}_K}\int_{0}^{||x||^{-1}_L}
\frac{r^{n-2}}{(1+\delta\ r^2)^{n-1}}drdx.
\end{eqnarray*}
But this follows from the assumption of the theorem, the Parseval's formula on the sphere
(Corollary \ref{Cor:Parseval}) and formula for the volume of central sections (Lemma
\ref{Lem:Main}). Indeed, let $\gamma_0$ be  the measure from Corollary \ref{Cor:Parseval}
corresponding to the Fourier transform of the positive definite distribution
$\frac{||x||^{-1}_K}{1+\delta\ (\frac{|x|}{||x||_K})^{2}}$, then
\begin{eqnarray*}
&&(2\pi)^n\int_{S^{n-1}}\frac{||x||^{-1}_K}{1+\delta\ ||x||^{-2}_K}\int_{0}^{||x||^{-1}_K} \frac{r^{n-2}}{(1+\delta\ r^2)^{n-1}}dr\, dx=\hspace{5cm}\\
&&=\int_{S^{n-1}}\left(\frac{||x||^{-1}_K}{1+\delta\ (\frac{|x|}{||x||_K})^{2}}\right)\cdot\left(
|x|^{-n+1}\int_{0}^{\frac{|x|}{||x||^{}_K}} \frac{r^{n-2}}{(1+\delta\ r^2)^{n-1}}dr\right) \, dx=\\
&&=\int_{S^{n-1}}\left(|x|^{-n+1}\int_{0}^{\frac{|x|}{||x||^{}_K}} \frac{r^{n-2}}{(1+\delta\ r^2)^{n-1}}dr\right)^\wedge(\theta) \, d\gamma_0(\theta)= \\
&&=\int_{S^{n-1}}\frac{\pi}{2^{n-1}} S_K(\theta)\, d\gamma_0(\theta)\le\int_{S^{n-1}}\frac{\pi}{2^{n-1}} S_L(\theta)\, d\gamma_0(\theta)= \\
&&=\int_{S^{n-1}}\left(|x|^{-n+1}\!\!\int_{0}^{\frac{|x|}{||x||^{}_L}}
\frac{r^{n-2}}{(1+\delta\ r^2)^{n-1}}dr\right)^\wedge  (\theta) \, d\gamma_0(\theta)=\\
&&= (2\pi)^n\int_{S^{n-1}}\frac{||x||^{-1}_K}{1+\delta\ ||x||^{-2}_K}\int_{0}^{||x||^{-1}_L}
\frac{r^{n-2}}{(1+\delta\ r^2)^{n-1}}drdx.
\end{eqnarray*}

 \qed

 \begin{Rmk}
Since $||x||_K^{-1}$ is positive definite for any convex origin-symmetric body in $\mathbb{R}^n$,
$n\le 4$ (see \cite{GKS}), the previous theorem implies the affirmative part of the original
Busemann-Petty problem in $\mathbb{R}^n$.
 \end{Rmk}

Now we investigate for which classes of bodies $\displaystyle\frac{||x||^{-1}_K}{1+\delta\
(\frac{|x|}{||x||_K})^2}$ is a positive definite distribution.

\begin{Prop}
Let $K$ be an origin-symmetric body in $B^n$, $n\le 4$.

i) If $K$ is h-convex then $\displaystyle\frac{||x||^{-1}_K}{1+ (\frac{|x|}{||x||_K})^2}$ is
positive definite.

ii) If $K$ is s-convex then $\displaystyle\frac{||x||^{-1}_K}{1- (\frac{|x|}{||x||_K})^2}$ is
positive definite.
\end{Prop}
\noindent{\bf Proof.} i) Consider a h-convex origin-symmetric body $K\subset B^n$, $n\le 4$. Define
a body $M$ by the formula:
$$||x||_M^{-1}=\frac{||x||^{-1}_K}{1+(\frac{|x|}{||x||_K})^2}.$$
It is enough to show that $M$ is e-convex. If we pass to polar coordinates then the map
$$\displaystyle (r,\theta)\mapsto \left(\frac{r}{1+r^2},\theta\right)$$ transforms the body $K$
into the body $M$.

Take two points in $K$ and connect them by a hyperbolic segment. This segment belongs to $K$ since
$K$ is h-convex. Consider the 2-dimensional plane through the origin and these 2 points. The
section of the body $K$ by this plane is a 2-dimensional h-convex body. Introduce polar coordinates
on this plane and (without loss of generality) assume that the h-geodesic segment has the equation
$r^2-a\ r \cos\phi +1=0$. Applying the above transformation one can see that this h-segment gets
mapped into an e-segment given by the equation $\displaystyle r=\frac{1}{a\cos\phi}$. Therefore the
body $M$ is e-convex and $(||x||_M^{-1})^\wedge$ is positive in dimensions $n\le 4$ (see
\cite{GKS}).

ii) Similar to (i). Take a s-geodesic given by the equation $r^2+a\ r \cos\phi - 1=0$. The image of
this geodesic under the map

\begin{eqnarray}\label{eqn:map}
\displaystyle (r,\theta)\mapsto\left(\frac{r}{1-r^2},\theta\right)
\end{eqnarray} is an e-geodesic
$\displaystyle r=\frac{1}{a\cos\phi}$.

\qed

Since every s-convex body containing the origin is h-convex, we have the following

\begin{Cor}
 $\displaystyle\frac{||x||^{-1}_K}{1+
(\frac{|x|}{||x||_K})^2}$ is positive-definite for every origin-symmetric s-convex body $K$ in dimension $n\le 4$.

This fact combined  with Theorem \ref{Thm:Main} implies the affirmative answer to the spherical
Busemann-Petty problem for $n\le 4$.
\end{Cor}

However not every h-convex body is s-convex and this idea will be used in constructing
counterexamples to the hyperbolic Busemann-Petty problem.

First we remind the following fact:

\begin{Thm}\label{Thm:GKS}{\rm (\cite{GKS}, Theorem 1)}
Let $K$ be an origin-symmetric star body in $\mathbb{R}^n$ with $C^\infty$ boundary, and let
$k\in\mathbb{N}\setminus \{0\}$, $k\ne n-1$. Suppose that $\xi\in S^{n-1}$, and let $A_\xi$ be the
corresponding parallel section function of $K$: $A_\xi(z)=\int_{K\cap \langle x,\xi\rangle=z} dx$.
\newline
(a) If $k$ is even, then
$$(||x||^{-n+k+1})^\wedge(\xi)=(-1)^{k/2}\pi (n-k-1) A_\xi^{(k)}(0).$$
\newline
(b) If $k$ is odd, then
\begin{eqnarray*}
(||x||^{-n+k+1})^\wedge(\xi) =(-1)^{(k+1)/2}2(n-1-k)k!\times\hspace{3.5cm}\\
\times\int_0^\infty
\frac{A_\xi(z)-A_\xi(0)-{A''}_\xi(0)\frac{z^2}{2}-\cdots-A_\xi^{(k-1)}(0)\frac{z^{k-1}}{(k-1)!}}{z^{k+1}}dz,
\end{eqnarray*}
where $A_\xi^{(k)}$ stands for the derivative of the order $k$ and the Fourier transform is
considered in the sense of distributions.
\end{Thm}
Now we can prove the following
\begin{Prop}
There exist h-convex origin-symmetric bodies in $B^n$, $n\ge 3$ that give a counterexample to the
hyperbolic  Busemann-Petty problem.

\end{Prop}
\noindent {\bf Proof.} In view of Theorem \ref{Thm:negative} we are interested only in the cases
$n=3$ and $4$.  First we construct a body $L$ for which
$\displaystyle\frac{||x||^{-1}_L}{1-(\frac{|x|}{||x||_L})^2}$ is not positive definite.

Let $L$ be a circular cylinder of radius $\sqrt{2}/2$ with  $x_1$ being its  axis of revolution.
(See Fig.2) To the top and bottom of  the cylinder attach spherical caps, that are totally geodesic
in the spherical metric. Clearly the body $L$ constructed this way is e-convex and therefore
h-convex.
% as a body in $\mathbb{H}^n$ (in the sense that it has nonzero principal curvatures of the same sign).
Using the formula
\begin{equation}\label{eqn:M}
||x||^{-1}_M=\displaystyle\frac{||x||^{-1}_L}{1-(\frac{|x|}{||x||_L})^2}
\end{equation} we define a body $M$.

%\centering

\begin{center}
\includegraphics{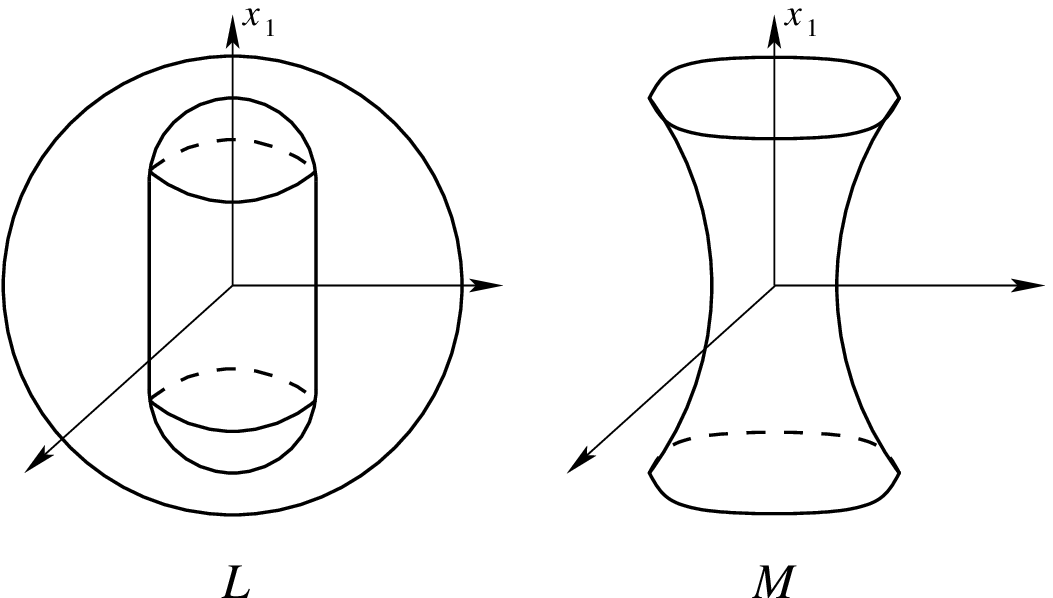}
%\nopagebreak

 Figure 2
\end{center}

Clearly the body $M$ is the image of $L$ under the map (\ref{eqn:map}). It can be checked directly
that the cylinder is mapped into the surface of revolution obtained by rotating the hyperbola
$x_2=\displaystyle \frac{1}{2}\left(\sqrt{2}+\sqrt{2+4x_1^2}\right)$ about the $x_1$-axis, and the
top and bottom spherical caps are mapped into flat disks.

In fact the body $L$ constructed above is not smooth. But we can approximate it by infinitely
smooth e-convex bodies that differ from $L$ only in a small neighborhood of the edges. Since the
body $M$ is obtained from  $L$ by (\ref{eqn:M}), and the denominator in (\ref{eqn:M}) is never
equal to zero, the body $M$ is also infinitely smooth. (Now that the bodies $L$ and $M$ are
smooth, Figure 2 might be confusing, but we wanted to make it as simple as possible, just to
emphasize the idea).

Now that we defined the body $M$, we can explicitly compute its parallel section function
$A_{M,\xi}$ in the direction of the $x_1$-axis.
\begin{eqnarray*}
A_{M,\xi}(t)=\left\{
\begin{array}{ll}
\displaystyle  \pi\left(\frac{\sqrt{2}+\sqrt{2+4t^2}}{2}\right)^2, & \mbox{ in dimension }n=3,\\
\displaystyle \frac{4\pi}{3} \left(\frac{\sqrt{2}+\sqrt{2+4t^2}}{2}\right)^3, & \mbox{ in dimension
}n=4.
\end{array}\right.
\end{eqnarray*}

Since $M$ is an infinitely smooth body,  $(||x||_M^{-1})^\wedge$ is a function. Applying Theorem
\ref{Thm:GKS} with $n=3$ and $q=1$ we get
$$(||x||_M^{-1})^\wedge(\xi)=-2 \int_0^\infty \frac{A_{M,\xi}(t)-A_{M,\xi}(0)}{t^2}dt.$$

Let the height of the cylindrical part of $L$ be equal to $\sqrt{2}-2\epsilon$ and the hight of its
image under (\ref{eqn:map}) equal to $N$. If $\epsilon$ tends to zero, the top and bottom parts of
the body $L$ get closer to the sphere $x_1^2+\cdots+x_n^2=1$. Recalling the definition of the
radial function of $M$:
$$\rho_M(x)=\displaystyle\frac{\rho_L(x)}{1-\rho_L(x)^2},\quad \forall x\in S^{n-1},$$
one can see that the the body $M$ becomes larger in the direction of $x_1$ as $\epsilon\to 0$, and
therefore its height $N$  approaches infinity.

Since in dimension $n=3$ the section function can be written as
$A_{M,\xi}(t)=\pi\left(1+t^2+\sqrt{1+2t^2}\right)$ for $-N\le t \le N$, we get:
\begin{eqnarray*}
(||x||_M^{-1})^\wedge(\xi)  &=&-2 \pi\int_0^N \frac{1+t^2+\sqrt{1+2t^2}-2}{t^2}dt-2
                                \pi\int_N^\infty\frac{(-2)}{t^2}dt\le\\
                            &\le&-2 \pi\int_0^N dt+4\pi\int_N^\infty\frac{1}{t^2}dt=\\
                            &=&-2\pi N+\frac{4\pi}{N}<0
\end{eqnarray*}
for $N$ large enough.

If $n=4$ and $q=2$ Theorem \ref{Thm:GKS} implies $$(||x||_M^{-1})^\wedge(\xi)=-\pi
A_{M,\xi}^{''}(0)< 0,$$ since the second derivative of the function $A_{M,\xi}$ in dimension $n=4$
equals: $A_{M,\xi}^{''}(0)=8\sqrt{2}\cdot\pi $.

 Thus  we have proved that
$\displaystyle\left(\frac{||x||^{-1}_L}{1-(\frac{|x|}{||x||_L})^2}\right)^\wedge(\xi)=(||x||_M^{-1})^\wedge(\xi)$
is negative for some direction $\xi$.

Now apply a standard argument to construct another body $K$ which along with the body $K$ provides
a counterexample to the hyperbolic Busemann-Petty problem (cf.  \cite{K3}, Theorem 2 or \cite{Zv},
Theorem 2). By continuity of $(||x||_M^{-1})^\wedge$ there is a neighborhood of $\xi$ where this
function is negative. Let
$$\Omega=\{\theta \in S^{n-1}: (||x||_M^{-1})^\wedge(\theta)<0\}.$$
Choose a non-positive infinitely-smooth even function $v$ supported on $\Omega$. Extend $v$ to a
homogeneous function $r^{-1}v(\theta)$ of degree $-1$ on $\mathbb{R}^n$. By Lemma 5 from \cite{K3}
we know that the Fourier transform of $ r^{-1}v(\theta)$ is equal to $r^{-n+1}g(\theta)$ for some
infinitely smooth function $g$ on $S^{n-1}$.

To construct a counterexample to the Busemann-Petty problem, define another body $K$ as follows:

\begin{eqnarray*}
\int_{0}^{||\theta||^{-1}_K} \frac{r^{n-2}}{(1-r^2)^{n-1}}dr=\int_{0}^{||\theta||^{-1}_L}
\frac{r^{n-2}}{(1-r^2)^{n-1}}dr+\epsilon g(\theta)
\end{eqnarray*}
for some  $\epsilon>0$  small enough (to guarantee that $K$ is still convex in hyperbolic sense).
Indeed, define a function $\alpha_\epsilon (\theta)$ such that
\begin{eqnarray*}
\int_{0}^{||\theta||^{-1}_L} \frac{r^{n-2}}{(1-r^2)^{n-1}}dr+\epsilon v(\theta)=
\int_{0}^{||\theta||^{-1}_L+\alpha_\epsilon (\theta)} \frac{r^{n-2}}{(1-r^2)^{n-1}}dr,
\end{eqnarray*}
then
\begin{eqnarray*}
||\theta||^{-1}_K= ||\theta||^{-1}_L+\alpha_\epsilon (\theta).
\end{eqnarray*}
Note that in our construction $L$ is e-convex, but we can perturb it a little (by adding $\alpha
|\theta|_2$  to the norm $||\theta||_L$ with $\alpha>0$ small enough), so  we can assume that L is
strictly e-convex. Therefore one can choose $\epsilon$ small enough such that $K$ is also e-convex
(for details see \cite{Zv}, Proposition 2). Hence we can assume that both $L$ and $K$ are h-convex.

Using Lemma \ref{Lem:Main}  we get
\begin{eqnarray*}
\mathrm{vol}_{n-1}(K\cap \xi^\perp)&=&\frac{2^{n-1}}{\pi}\left(|x|^{-n+1}\int_{0}^{|x|/||x||^{}_K}
\frac{r^{n-2}}{(1-r^2)^{n-1}}dr\right)^\wedge(\xi)=\\
&=&\frac{2^{n-1}}{\pi}\left(|x|^{-n+1}\int_{0}^{|x|/||x||^{}_L}
\frac{r^{n-2}}{(1-r^2)^{n-1}}dr\right)^\wedge(\xi)+ \epsilon
 v(\xi)\le \\
&\le&\frac{2^{n-1}}{\pi}\left(|x|^{-n+1}\int_{0}^{|x|/||x||^{}_L}
\frac{r^{n-2}}{(1-r^2)^{n-1}}dr\right)^\wedge(\xi)= \\
&=&\mathrm{vol}_{n-1}(L\cap \xi^\perp).
\end{eqnarray*}

Proceeding as in the proof of Theorem \ref{Thm:Main} we can show the opposite inequality for
volumes. Since the body $L$ is infinitely smooth, one can use the Parseval's formula in the form of
Lemma \ref{Lem:Parseval}:
\begin{eqnarray*}
&&(2\pi)^n\int_{S^{n-1}}\frac{||x||^{-1}_L}{1-\ ||x||^{-2}_L}\int_{0}^{||x||^{-1}_K} \frac{r^{n-2}}{(1-\ r^2)^{n-1}}dr dx=\\
&&=\int_{S^{n-1}}\left(\frac{||x||^{-1}_L}{1-\ (\frac{|x|}{||x||_L})^{2}}\right)^\wedge(\theta)
\left(|x|^{-n+1}\int_{0}^{\frac{|x|}{||x||^{}_K}} \frac{r^{n-2}}{(1-\ r^2)^{n-1}}dr\right)^\wedge(\theta) d\theta=\\
&&=\int_{S^{n-1}}\left(\frac{||x||^{-1}_L}{1-\ (\frac{|x|}{||x||_L})^{2}}\right)^\wedge(\theta)
\left(|x|^{-n+1}\int_{0}^{\frac{|x|}{||x||^{}_L}}
\frac{r^{n-2}}{(1-\ r^2)^{n-1}}dr\right)^\wedge(\theta) d\theta+\\
&&+  \int_{S^{n-1}}\left(\frac{||x||^{-1}_L}{1-\
(\frac{|x|}{||x||_L})^{2}}\right)^\wedge(\theta)\cdot
\epsilon v\left(\theta\right)d\theta>\\
&& >(2\pi)^n\int_{S^{n-1}}\frac{||x||^{-1}_L}{1-\ ||x||^{-2}_L}\int_{0}^{||x||^{-1}_L}
\frac{r^{n-2}}{(1-\ r^2)^{n-1}}dr dx.
\end{eqnarray*}

\qed

{\bf Acknowledgments}. The author wishes to thank A.Koldobsky for useful discussions and A.Zvavitch
for many suggestions and bringing to my attention his results from \cite{Zv}.

\end{document}